\documentclass{amsart}
\usepackage{amsmath,amsthm,amssymb,amsfonts}
\usepackage[colorlinks=true]{hyperref}

\hypersetup{urlcolor=blue, citecolor=blue}

\def\arXiv#1{   {\href{http://arxiv.org/abs/#1}
   {{\mdseries\ttfamily arXiv:#1}}}}\let\arxiv\arXiv

\let\MR\mr

\def\doi#1{   {\href{http://dx.doi.org/#1}
   {{\mdseries\ttfamily DOI}}}}

\newcommand{\al}{\alpha}    \newcommand{\be}{\beta}
\newcommand{\de}{\delta}    
\newcommand{\vep}{\epsilon}  \newcommand{\ep}{\varepsilon}
   \newcommand{\La}{\Lambda}

\newcommand{\ga}{\gamma}    
\newcommand{\R}{\mathbb{R}}\newcommand{\Z}{\mathbb{Z}}
\newcommand{\N}{\mathbb{N}}

\newcommand{\pt}{\partial_t}\newcommand{\pa}{\partial}
\newcommand{\les}{{\lesssim}}
\newcommand{\beeq}{\begin{equation}}\newcommand{\eneq}{\end{equation}}

\newcommand{\ang}{{\not\negmedspace\nabla}}

\newcommand{\Sp}{{\mathbb S}}

\newenvironment{prf}{\noindent {\bf Proof.} }{\endprf\par}
\def \endprf{\hfill  {\vrule height6pt width6pt depth0pt}\medskip}

\numberwithin{equation}{section}

\newcommand{\gm}{\mathfrak{g}} 
 \newcommand{\eps}{\varepsilon}

\def\<{\langle}             \def\>{\rangle}
\def\({\left(}                 \def\){\right)}

\newtheorem{thm}{Theorem}[section]

\newtheorem{lem}[thm]{Lemma}

\theoremstyle{remark}
\newtheorem{rem}{Remark}

\theoremstyle{definition}

\title[Glassey conjecture on manifolds]
{The Glassey conjecture on asymptotically flat manifolds}

\author{Chengbo Wang}
\address{Department of Mathematics\\
                Zhejiang University\\
                Hangzhou 310027, China}
\email{wangcbo@gmail.com}
\urladdr{http://www.math.zju.edu.cn/wang}
\thanks{
The author was supported by Zhejiang Provincial Natural Science Foundation of
China LR12A01002, the Fundamental Research Funds for the Central Universities, NSFC 11301478, 11271322 and J1210038.
}
\date{}
\dedicatory{} \commby{}

\begin{document}

\begin{abstract}
We verify the $3$-dimensional Glassey conjecture on asymptotically flat manifolds $(\R^{1+3},\gm)$, where the metric $\gm$ is certain small space-time perturbation of the flat metric, as well as the nontrapping asymptotically Euclidean manifolds.
Moreover, for radial asymptotically flat manifolds $(\R^{1+n},\gm)$ with $n\ge 3$, we verify the Glassey conjecture in the radial case. High dimensional wave equations with higher regularity are also discussed. The main idea is to exploit local energy and KSS estimates with variable coefficients, together with the weighted Sobolev estimates including trace estimates.
\end{abstract}



\maketitle 

\section{Introduction}
\label{sec-Intro}

The purpose of this article is to study the long time existence of solutions for the Cauchy problem of the nonlinear wave equations of the type
$$\Box u=c_1 |u_t|^p+c_2 |\nabla_x u|^p
$$ on certain asymptotically flat manifolds, with small initial data.

In the 1980's, Glassey made a conjecture that the critical exponent $p$ for the problem, posed on the Minkowski space-time, to admit global solutions with small, smooth initial data with compact support is $$p_c=1+\frac{2}{n-1}$$ in \cite{Glassey} (see also Schaeffer \cite{Sch86}, Rammaha \cite{Ram87}), where $n$ is the spatial dimension.
The conjecture was verified for dimension $n=2, 3$ for general data (Hidano and Tsutaya \cite{HiTs95} and Tzvetkov \cite{Tz98} independently, as well as the radial case in Sideris \cite{Si83} for $n=3$). For higher dimension $n\ge 4$, when $c_2=0$ and $c_1>0$, the blow up results (together with an explicit upper bound of the lifespan) for $p\le p_c$ were obtained in Zhou \cite{Zh01} (see also Zhou and Han \cite{ZhHa10} for problems with time-independent compact metric perturbation). Recently, for the radial data, the existence results with sharp lifespan for any $p\in (1, 1+2/(n-2))$ was proved in Hidano, Yokoyama and  the author \cite{HWY11} (see also Fang and the author \cite{FaWa10} for the critical case $n=2$ and $p=3$), which particularly proved the Glassey conjecture in the radial case.

In this paper, we prove analogs of the results on the Glassey conjecture
of \cite{HiTs95}, \cite{Tz98} and \cite{HWY11}, on certain asymptotically flat manifolds.
One of the main ingredients in the proof is the local energy estimates with variable coefficients of Metcalfe and Sogge \cite{MeSo06} and Hidano, Yokoyama and the author \cite{HWY10}. The local energy estimates first appeared in Morawetz \cite{Mo2} and subsequently in many works, see e.g. \cite{Strauss75, KPV, SmSo, KSS02, Sterb, HY05, MeSo06} and references therein. Such estimates are known to be fairly robust and variants for metric perturbations
were proved in e.g. \cite{Burq, Al06, MeSo06, MeTa07, BoHa, SoWa10, HWY10, HWY11, LMSTW} and references therein.
Another key ingredient is the weighted Sobolev estimates, including the trace estimates. Together with the local energy estimates, such kind of estimates have been proved to be very useful,
see e.g. \cite{Kl1, KSS02,  MeSo06, HMSSZ, FaWa11, LMSTW}.

Let us begin with the space-time manifolds we will work on. As usual, we use $(x^0,x^1,\cdots ,x^n)=(t,x)$ to denote points in $\R^{1+n}$, and $\partial_\al={\partial}/{\partial x^\al}$ as partial derivatives,
with the abbreviations $\partial=(\partial_0,\partial_1,\cdots ,\partial_n)=(\partial_t,\nabla)$.
We consider the asymptotically flat Lorentzian manifolds $(\R^{1+n} , \gm)$ with
$$\gm =  g_{\al\be} (t,x) \, d x^\al \, d x^\be =\sum_{\al,\be=0}^n
 g_{\al\be} (t,x) \, d x^\al \, d x^\be .$$
Here, we have used the convention that Greek indices $\al$, $\be$, $\ga$ range from
$0$ to $n$ and Latin indices $i$, $j$, $k$ from $1$ to $n$. We will also use $a$, $b$, $c$ to
denote multi-indices. Moreover, the Einstein summation convention will be performed over repeated indices.

We will consider two types of manifolds. The first one is small, time-dependent, asymptotically flat perturbation of the Minkowski metric. More precisely, with Euclidean radius $r=\sqrt{|x|^2}$, we assume $g_{\al\be} \in C^{\infty} ( \R^{1+n} )$ and, for some fixed $\rho_1>0$, $\rho_2>1$ and $\delta\ll 1$,
\begin{equation}\tag{H1} \label{H1}
\gm=m+g_1(t,r)+g_2(t,x),
\end{equation}
where $(m_{\al\be})=Diag(1,-1,-1,\cdots,-1)$ is the standard Minkowski metric, $\<x\>=\sqrt{1+|x|^2}$,
\beeq\tag{H1.1} \label{H1.1}|\partial^a_{t,x} g_{i,\al\be} |\les_a \delta  \< x \>^{- \vert a \vert - \rho_i}, i=1,2, \rho=\min(\rho_1,\rho_2-1), \rho_1<\rho_2,\eneq
and moreover, we assume the first perturbation $g_1$ is radial. By radial metric, we mean that, when writing out the metric $\gm$, with $g_2=0$, in polar coordinates $(t,x)=(t,r\omega)$ with $\omega\in\mathbb{S}^{n-1}$, we have
$$m+g_1(t,r)=\tilde g_{00}(t,r)dt^2+2\tilde g_{01}(t,r)dtdr+\tilde g_{11}(t,r)dr^2+\tilde g_{22}(t,r) r^2 d\omega^2\ .$$
In this form, the assumption \eqref{H1.1} for $g_1$ on asymptotically flatness is equivalent to the following requirement
\begin{equation}\tag{H1.2} \label{H1.2}
|\partial^a_{t,x} ( \tilde g_{00}-1,  \tilde g_{11}+1, \tilde  g_{22}+1, \tilde  g_{01})|
\les_a \delta  \< x \>^{- \vert a \vert - \rho}  .
\end{equation}
 An example of such metric can be $$g_{\al\be}=m_{\al\be}+\de \<x\>^{-\rho}+\de\phi(t/\<x\>) \<x\>^{-\rho}$$ with $\phi\in C_0^\infty(\R)$.
 Here, and in what follows, we use $A\les B$ to stand for $A\le C B$ where the constant $C$ may change from line to line.
 When $\delta\ge 0$ is small enough, it is clear that the metric $\gm$ is a nontrapping perturbation. 
Notice that the form of the metrics mimic to that in Tataru \cite{Ta13} and Metcalfe, Tataru and Tohaneanu \cite{Ta12}, and $g_1$ could be long range perturbation. 
Also, the metric models the black hole metrics occurring in mathematical relativity,  near spatial infinity. For example, for large enough radius $r$, Schwarzschild metric is radial metric \eqref{H1} with $\rho_1=1$, and Kerr metric is a metric \eqref{H1} with $\rho_1=1$ and $\rho_2=2$. The separation of the radial part from the general non-radial metric is crucial in many problems, see e.g. \cite{Ta13, Ta12, LMSTW}.

We will also consider nontrapping asymptotically Euclidean manifolds. That is, we assume that
\begin{equation}\tag{H2} \label{H2}
\gm=m+g_1(r)+g_2( x), \gm \text{ is nontrapping},
\end{equation}
where
we suppose $g_1$ and $g_2$ are of the form $g_{jk}dx^jdx^k$,
\begin{equation}\tag{H2.1} \label{H2.1}
|\nabla^a_{x} g_{i,jk} |\les_a   \< x \>^{- \vert a \vert - \rho_i}, i=1,2, \rho=\min(\rho_1,\rho_2-1), \rho_1<\rho_2
\end{equation}
and as \eqref{H1.2}, we assume the first perturbation $g_1$ is radial.

 With $|\gm| =(-1)^n \det (g_{\al\be})$, the d'Alembertian operator associated with $\gm$ is given by
$$\Box_{\gm} =  \sqrt{|\gm|}^{-1} \partial_\al g^{\al\be} \sqrt{|\gm|} \partial_\be ,
$$
where $(g^{\al\be} (t,x))$ denotes the inverse matrix of $(g_{\al\be}(t,x))$.

Let $n\ge 3$, $p>1$, we consider initial value problems for the following nonlinear wave equations, 
\begin{equation}
\label{a-pde}
\left\{\begin{array}{l}\Box_\gm u = \sum_{\al=0}^n a_\al(u) |\pa_\al u|^p\equiv F_p(u) ,\ x\in\R^n\\
u(0,x)=f(x), \pt u(0,x)=g(x)\ ,
\end{array}\right.
\end{equation}
for given smooth functions $a_\al$,
as well as the radial problems (with $g_2=0$)
\begin{equation}
\label{a-pde2}
\left\{\begin{array}{l}\Box_\gm u = c_1 |\pt u|^p+c_2 |\nabla_x u|^p\equiv G_p(u)\ , x\in\R^n \\
u(0,x)=f(x), \pt u(0,x)=g(x)\ ,
\end{array}\right.
\end{equation} for constants $c_1$, $c_2$.
We will study the long time existence of such problems with small enough initial data, according to certain norm.

Before stating our results, let us give some more notations.
At first, the vector fields to be used will be labeled as
\[
Y=(Y_1,\cdots,Y_{n(n+1)/2})=(\nabla,\Omega)
\]
with rotational vector fields
$$\Omega_{ij}=x_i \pa_j-x_j\pa_i, 1\le i<j\le n\ . $$
For fixed $T>0$, the space-time norm $L^q_T L^r_x$ is simply $L^q_t([0,T], L^r_x (\R^n))$. In the case of $T=\infty$, we use $L^q_t L^r_x$ to denote $L^q_t([0,\infty), L^r_x (\R^n))$. As usual, we use
$\|\cdot\|_{E_m}$  to denote the energy norm of order $m\ge 0$,
\beeq\label{eq-E}\|u\|_{E}=\|u\|_{E_0}=\|\pa u\|_{L^\infty_T L^2_x}\ ,
\|u\|_{E_{m}}=\sum_{|a|\le m}\| Y^a u\|_{E}\ .\eneq
We will use $\|\cdot\|_{LE}$ to denote the (strong) local energy norm
\beeq\label{eq-LE}
\|u\|_{LE}=\|u\|_E+\|\pa u\|_{l^{-1/2}_\infty (L^2_t L^2_x)}+\|r^{-1} u\|_{l^{-1/2}_\infty (L^2_t L^2_x)},
\eneq
 where  we write
 \[
 \|u\|_{l^s_q(A)} = \|(\phi_j(x)u(t, x))\|_{l^s_q(A)},
\]
for a partition of unity subordinate to the dyadic (spatial) annuli,
$\sum_{j\ge 0}\phi^2_j(x)=1$.
The local energy includes the energy norm and could control the KSS norm, see Lemma \ref{thm-KSS}.
On the basis of the local energy norm, we can similarly define
$\|u\|_{LE_m}$,  and the dual norm $LE^*=l^{1/2}_1 L^2_t L^2_x$.

In proving the radial Glassey conjecture, we need to obtain stronger control on the local information, in contrast to the local energy norm \eqref{eq-LE}. For this purpose, we introduce a modified local energy norm, mimicking that occurred in Hidano and Yokoyama \cite{HY05} (see also \cite{HWY10}),
\beeq\label{eq-LE2}\|u\|_{\widetilde{LE}}=\left\|r^{-1/2+\mu}\<r\>^{-\mu'} \pa u\right\|_{L^2_t L^2_x}+\left\|r^{-3/2+\mu}\<r\>^{-\mu'} u\right\|_{L^2_t L^2_x}+\|\pa u\|_{L^\infty_t L^2_x}\eneq
with $\mu\in (0,1/2)$ and $\mu'>\mu$ to be specified.
On the basis of $\widetilde{LE}$, we can similarly define
$$\|u\|_{\widetilde{LE}_{m}}=\sum_{|a|\le m}\| \nabla^a u\|_{E},\ 
\|F\|_{\widetilde{LE}^*}=\|r^{1/2-\mu} \<r\>^{\mu'} F\|_{L^2_t L^2_{x}}\ .
$$
Here, we do not use the vector field $\Omega$, since we will only apply such norms in the radial problems.

We can now state our main results. The first result is about the problem \eqref{a-pde} with general data,
which verifies the $3$-dimensional Glassey conjecture on asymptotically flat manifolds and nontrapping asymptotically Euclidean manifolds.
\begin{thm}
\label{thm-3}
Consider the problem \eqref{a-pde} on the manifold $( \R^{1+3} , \gm)$ satisfying \eqref{H1} or \eqref{H2} with $\rho=\min(\rho_1,\rho_2-1)>0$. Let $p>2$,
there exist small positive constants $\eps_0$ and $\de_0$, such that the Cauchy problem with $\de\le \de_0$ has a unique global solution $u\in C([0,\infty); H^3(\R^3))\cap C^1([0,\infty); H^2(\R^3))$ when the initial data
satisfy
\begin{equation}
\label{a-smallness}
\sum_{|a|\le 2}\|\pa Y^a u(0)\|_{L^2(\R^3)}=\ep \le  \ep_0, \ \|u(0)\|_{L^2(\R^3)}<\infty\ .
\end{equation} Moreover, it satisfies
$\|u\|_{LE_{2}}\les \ep$.
\end{thm}

As a byproduct, provided that the nonlinearity is sufficient smooth, we can also deal with the case $p\ge 2$ and $n\ge 3$, which generalizes the works for $p=2$ on asymptotically Euclidean manifolds \eqref{H2} (with $g_1=0$) in Bony and H\"{a}fner \cite{BoHa}, Sogge and the author \cite{SoWa10}.
For small asymptotically flat manifolds \eqref{H1} (with $g_1=0$) and $p=2$,
the almost global existence and global existence for the solutions was essentially proved even for the quasilinear problems, in Yu and the author \cite{WY12}.
 Notice also that we have improved the condition on the regularity. For simplicity, we only consider the case of constant coefficients ($a_\al(u)=a_\al$).
\begin{thm}
\label{thm-1}
Let $n\ge 3$,  $p\ge 2$, $m=[\frac{n+2}2]$, $p> m$ or $p\in \N$ and $F_p(u)$ smooth in $\pa u$ (e.g. $p=2$ and $F_p(u)=|\pa u|^2$, $p=3$ and $F_p(u)=|\pa u|^2 \pt u$). Consider the problem \eqref{a-pde} with $a_\al(u)=a_\al$ on the manifold $(\R^{1+n} , \gm)$ satisfying \eqref{H1} or \eqref{H2} with $\rho=\min(\rho_1,\rho_2-1)>0$.
Then there exist small positive constants $\eps_0$ and $\de_0$, such that for any data
 \begin{equation}
\label{a-smallness2}
\sum_{|a|\le m}\|\pa Y^a u(0)\|_{L^2(\R^n)}=\ep \le  \ep_0, \ \|u(0)\|_{L^2(\R^n)}<\infty\ ,
\end{equation}
the Cauchy problem with $\de\le \de_0$ has a unique global solution $u\in C([0,T_*); H^{m+1}(\R^n))\cap C^1([0,T_*); H^m(\R^n))$ with $\|u\|_{LE_{m}}\les \ep$. Here $T_*=\infty$, except for the case $n=3$ and $p=2$, where $T_*=\exp(c/\ep)$ for some constant $c>0$.
\end{thm}

Turning to the problem \eqref{a-pde2} with radial data, we can prove the long time existence of the radial solutions, in spirit of \cite{HWY11}. To simplify the exposition, we will only prove a weaker version of the global existence theorem, comparing that in \cite{HWY11}. It is not hard to see that our proof could be 
adapted to prove the same set of Theorems 1.1-1.3 (with $n\ge 3$) in \cite{HWY11}, by noticing Lemma \ref{thm-KSS-HY-2} and using $$\La_2=\|\pa \nabla u(0)\|_{L^2(\R^n)}+\de \|\pa u(0)\|_{L^2_x}\ .$$ We leave the details to the reader.

We shall use $H^m_{\rm{rad}}$ to denote the space of spherically symmetric functions in the usual Sobolev space $H^m$. We have the following global existence theorem for $p>p_c$ and $n\ge 3$, which can be viewed as a  positive solution for the radial Glassey conjecture in the setting of asymptotically flat manifolds.
\begin{thm}\label{thm-sub}
Let $n\ge 3$ and $1+2/(n-1)<p<1+2/(n-2)$.
Consider the problem \eqref{a-pde2} with initial data $(f, g)\in H^2_{{\rm{rad}}}\times H^1_{{\rm{rad}}}$, posed on the manifold $( \R^{1+n} , \gm)$ satisfying \eqref{H1} with $g_2=0$ and $\rho>0$.
Then there exist constants $\ep_0, \de_0>0$, such that if
$$\|\pa u(0)\|_{H^1}=\ep \le \ep_0\ ,\ \de\le\de_0\ , $$
then we have a unique global solution $u$ to \eqref{a-pde2} satisfying
$$ u\in C([0,\infty); H^2_{\rm{rad}}(\R^n))\cap C^1([0,\infty); H^1_{\rm{rad}}(\R^n))\ ,$$
 $$\|u\|_{\widetilde{LE}_1} \les \ep\ , \mu=\frac{1}2- \frac{n-2}4 (p-1)\ ,\ \mu'=\frac{p-1}4\ .$$
\end{thm}
In Theorem \ref{thm-sub}, the technical restriction $p<1+2/(n-2)$ is partly due to the regularity, since we are assuming $H^2$ regularity of the data, while the critical scaling regularity for the problem is $s_c=n/2+1-1/(p-1)$, and $s_c<2$ if and only if $p<1+2/(n-2)$.
We could expect to obtain more general results if we relax the regularity assumption. However, we could not expect to completely overcome this difficulty just by increasing the regularity, due to the limited regularity of the nonlinearity (noticing that $1+2/(n-2)<2$ for $n>4$ and so $|x|^{1+2/(n-2)}\not\in C^2$). Because of this, current technology can only be adapted for some specific situations. 
\begin{thm}
\label{thm-add}
Let $n\ge 4$ be even numbers and $p> n/2$.
Consider the problem \eqref{a-pde2}, posed on the manifold $( \R^{1+n} , \gm)$ satisfying \eqref{H1} with $g_2=0$ and $\rho>0$. When $(u(0),\pt u(0))\in H^{n/2+1}_{\rm{rad}}\times H^{n/2}_{\rm{rad}}$,
$$\|\pa u(0)\|_{H^{n/2}}=\ep\ll 1, \|\pa u(0)\|_{L^\infty_x}\les 1$$ and $\de\ll 1$,
we have global existence and uniqueness of the radial solutions for \eqref{a-pde2} in $C([0,\infty);H^{{n/2}+1}_{{\rm{rad}}})\cap C^1([0,\infty); H^{n/2}_{\rm{rad}})$ with the property
$$\sum_{|a|\le n/2}\|\pa^a u\|_{\widetilde{LE}} \les \ep\ , \mu=\frac{1}4\ ,\ \mu'=\frac{(n-1)(p-1)-1}4 .$$
In addition, for $n=5$ and $p\in [2,3)$, we have a similar result in $C([0,\infty);H^{3}_{{\rm{rad}}})\cap C^1([0,\infty); H^2_{\rm{rad}})$, for small $H_{{\rm{rad}}}^{3}\times H_{{\rm{rad}}}^2$ data (without the boundedness assumption on $\pa u(0)$),
 with $\mu= \frac{3-p}4$ and $\mu'=\frac{3}4(p-1)$.
\end{thm}
\begin{rem}
Comparing Theorems \ref{thm-sub}, \ref{thm-add} and \ref{thm-1} applied to radial solutions, we see that we could prove global existence of radial solutions for $p>1+2/(n-1)$, except the case $p\in [5/3,2)$ for $n=5$ and $p\in [1+2/(n-2),[(n+1)/2])$ for $n\ge 6$. 
As we have mentioned, it is mainly due to the limited regularity of the nonlinearity.
It will be interesting to fill this gap. By exploiting Strichartz estimates, one may shrink the gap further.  However, it seems that the known Strichartz estimates for $\Box_\gm$ have more restriction on the assumption of the metric than \eqref{H1} (for example, an additional requirement $|g_1(t,r)|\les \de r^\mu$ for some $\mu>0$ may be required, see e.g. Metcalfe and Tataru \cite{MeTa07} and references therein).
\end{rem}

\begin{rem}
One may ask if we can prove similar results as Theorem \ref{thm-sub} on nontrapping asymptotically Euclidean manifolds, as in \cite{BoHa}, \cite{SoWa10} and \cite{WaYu11} for other related problems. However, it seems to us that the current available local energy estimates are not enough to control all of the required local information.
\end{rem}

This paper is organized as follows. In the next section, we give several Sobolev type estimates, including the trace estimates.
In Section \ref{sec-KSS}, we give various forms of the local energy estimates and KSS estimates for $\Box_\gm$ satisfying \eqref{H1} or \eqref{H2}, together with the higher order estimates.
The remaining sections are devoted to the proof of the Theorems.

\section{Sobolev-type estimates}
\label{sec-Sobolev}
In this section, we give several Sobolev type estimates, including the trace
estimates. At first, let us recall the trace estimates (see (1.3), (1.7) in Fang and the author \cite{FaWa11} and references therein)
\begin{lem}[Trace estimates]\label{thm-trace}
Let $n\ge 2$ and $1/2< s<n/2$, then
\beeq\label{eq-trace}
\|r^{n/2-s} u\|_{L^\infty_r L^2_\omega}\les_s \|u\|_{\dot H^s}\ ,
\|r^{(n-1)/2} u\|_{L^\infty_r L^{2}_\omega}\les_s \|u\|_{ H^s}\ ,
\eneq
in particular, for radial functions, \beeq\label{eq-trace-rad}
\|r^{n/2-s}\<r\>^{s-1/2} u\|_{L^\infty_x}\les_s \|u\|_{ H^s_{\rm{rad}}}\ .
\eneq
\end{lem}

We will also need the following variant of the Sobolev embeddings.
\begin{lem}\label{thm-Sobo}
Let $n\ge 2$.
For any $m\in \R$ and $k\ge n/2-n/q$ with $q\in [2,\infty)$, we have
\beeq\label{eq-Sobo}
\|\<r\>^{(n-1)(1/2-1/q)+m} u\|_{L^q(\R^n)}\les \sum_{|a|\le k} \|\<r\>^{m} Y^a u\|_{L^2(\R^n)}\ .
\eneq
Moreover, we have
\beeq\label{eq-Sobo2}
\|\<r\>^{(n-1)/2+m} u\|_{L^\infty(\R^n)}\les \sum_{|a|\le [(n+2)/2]} \|\<r\>^{m}  Y^a u\|_{L^2(\R^n)}\ ,
\eneq
where $[k]$ stands for the integer part of $k$.
\end{lem}
\begin{prf}
We employ a similar proof of related estimates as in Lemma 3.1 of Lindblad, Metcalfe, Sogge, Tohaneanu and the author \cite{LMSTW}. For $q<\infty$,
by Sobolev's lemma $H^k\subset L^q$ on ${\mathbb R}\times \Sp^{n-1}$, we have for each $j\ge 0$ the uniform bounds
$$ \Bigl(\int_{j+1}^{j+2}\int_{\Sp^{n-1}}|v|^q \, d\omega dr\Bigr)^{\frac1q}
\lesssim \sum_{|a|\le k}\Bigl(\int_{j}^{j+3}\int_{\Sp^{n-1}}|Y^a v|^2 \, d\omega dr\Bigr)^{\frac12}.
$$ Hence,
$$ \|v\|_{L_x^{q}(r\in [j+1, j+2])} \lesssim j^{-(n-1)(1/2-1/q)} \sum_{|a|\leq k} \|
Y^a v\|_{L_x^{2}(r\in [j,j+3])} .$$ Or more generally, for any $m\in\R$,
\[\|r^{m+(n-1)(1/2-1/q)} v\|_{L^q_x(r\in [j+1,j+2])} \lesssim   \sum_{|a|\le k} \|r^{m} Y^a v\|_{L^2_x(r\in [j,j+3])}.\]
The factor $j^{-(n-1)(1/2-1/q)}$ on the right comes from the fact that the volume element for ${\mathbb R}^n$ is $r^{n-1} dr d\omega$. Recalling that the Sobolev's lemma on $\R^n$ also gives us
\[
\|\<r\>^{m+(n-1)(1/2-1/q)} v\|_{L^q_x(r\in [0,1])} \lesssim   \sum_{|a|\le k} \|\<r\>^{m} Y^a v\|_{L^2_x(r\in [0,2])}.
\]
The above two inequalities together imply \eqref{eq-Sobo} if we  $l^q$-sum over $j\ge 0$ using the
Minkowski integral inequality. A direct modification of the proof yields the case $q=\infty$, \eqref{eq-Sobo2},
which completes the proof.
\end{prf}

In the last section, we will also need the weighted Hardy-Littlewood-Sobolev inequalities (Stein and Weiss \cite{SW}). Here, we just record a particular case which we will use.
\begin{lem}\label{thm-wHLS}
Let $n\ge 2$.
For any $k\in (0,n)$ and $\al\in (-n/2, n/2-k)$,
we have
\beeq\label{eq-wHLS}
\|r^{-\al-k} u\|_{L^2(\R^n)}\les \|r^{-\al} \sqrt{-\Delta}^k u\|_{L^2(\R^n)}\ .
\eneq
In particular, when $k=2l$ with $l\in \Z_+$, we have
\beeq\label{eq-wHLS2}
\|r^{-\al-k} u\|_{L^2(\R^n)}\les \|r^{-\al} (-\Delta)^l u\|_{L^2(\R^n)}\ .
\eneq
\end{lem}

\section{Local energy estimates}\label{sec-KSS}
 In this section, we collect the required local energy estimates for $\Box_\gm$, together with the higher order estimates.

\subsection{Local energy estimates}

\subsubsection{Asymptotically flat manifolds}
\begin{lem}\label{thm-LE-g1-620}
Let $n\ge 3$ and consider the linear problem $\Box_\gm u=F$ on the manifold $(\R^{1+n} , \gm)$ satisfying \eqref{H1} with $\rho=\min(\rho_1,\rho_2-1)>0$. Then there exists a constant $\de_0$, such that for any $\de\le \de_0$, we have the following local energy estimates,
\beeq\label{eq-LE-g1-620}
\|u\|_{LE}\les \|\pa u(0)\|_{L^2_x}+\|F\|_{LE^*+L^1_t L^2_x}\ .
\eneq
\end{lem}

This result was essentially proved in Section 2 of  \cite{HWY10} and \cite{MeSo06}, see also Lemma 3.2 of \cite{HWY11} and \cite{WY12}. Here, let us record one version of the local energy estimates, which are consequences of the classical positive commutator method (also known as Friedrichs' abc method).
\begin{lem}\label{thm-KSS-0}
Let $f=f(r)$ be any fixed differential function.
For any solution $u\in C^\infty([0,T], C_0^\infty(\R^n))$ to the equation \beeq\label{612-1}(\pt^2-\Delta +h^{\al\be}(t,x)\pa_\al\pa_\be)u=F\eneq in $S_T=[0,T]\times \R^n$ with $h^{\al\be}=h^{\be\al}$, $\sum_{0\le \al,\be\le n}|h^{\al\be}|\le 1/2$ and $n\ge 3$,
we have
  \beeq\label{520-2}\int_{S_T} Q dxdt=\int_{S_T} f F \left(\pa_r+\frac{n-1}{2r}\right) u  dx dt-\int_{\R^n} P^0(t, \cdot)dx|_0^T\ ,
  \eneq
  where
$P^0=f  (m+h)^{0\beta}  \partial_\beta u \left(\pa_r+\frac{n-1}{2r}\right) u$, 
\beeq\label{520-4}Q_0=  \frac{ 2f-r f' }r \frac{|\ang u|^2}2+f' \frac{|\pa_r u|^2+
|\pt u|^2}{2}
-\frac{n-1}4\Delta \left(\frac{f }{r}\right) u^2\ ,\eneq
and
\begin{eqnarray*}
Q&=&Q_0+
\left(\frac{f}{2}\pa_r  h^{\alpha \beta} +\frac{f'}2 h^{\alpha \beta}
\right) \pa_\al u \pa_\be u
- f\pa_\ga h^{\ga\be} \pa_r u \pa_\be u
-
h^{i\be}  \frac{f}{r}
\not\negmedspace\pa_i u \pa_\be u
\\
&&
-
\omega_i h^{i\be}
f' \pa_r u \pa_\be u-\frac{n-1}2 \left(\omega_i h^{i\be}  \left(\frac{f}r\right)'
+\frac f r\pa_\ga h^{\ga \beta}   \right)
u \pa_\be u\nonumber
\end{eqnarray*}
with $\omega_i=x^i/r$, $\not\negmedspace\pa_i u=\pa_i u-\omega_i \pa_r u$, $|\ang u|^2=|\nabla u|^2-|\pa_r u|^2$.
\end{lem}
This is essentially coming from multiplying $f(r) \left(\pa_r+\frac{n-1}{2r}\right) u$ to the wave equation and a tedious calculation of integration by parts. See e.g. \cite{MeSo06}  P200 (5.4). Similar to the energy estimates, Lemma \ref{thm-KSS-0} is robust enough for us to obtain many interesting and useful estimates, particularly local energy and KSS estimates.

Typically, $f$ is chosen to be differential functions
satisfying
\beeq\label{520-1}f\le 1, 2f \ge r f'(r)\ge 0,  -\Delta (f/r)\ge 0,\eneq
which ensure that $Q_0$ is positive semidefinite.
In literatures, the typical choices including $f=1$ \cite{Mo2}, $1-(3+r)^{-\delta}$ ($\delta>0$) \cite{Strauss75}, $r/(R+r)$ \cite{Sterb, MeSo06}, $(r/(R+r))^{2\mu}$ ($\mu\in (0,1/2)$, \cite{HWY10, HWY11}).

With the help of Lemma \ref{thm-KSS-0} with $f_j=r/(2^j+r)$ and $j\ge 0$, and the classical energy estimates, we can give the proof of Lemma \ref{thm-LE-g1-620}.

At first, we notice that
\beeq\label{eq-620-1}\Box_\gm= \Box+ (g^{\al\be}-m^{\al\be})\pa_\al\pa_\be+r_1 \pa \ ,\eneq
here, and in what follows, we use $r_m$ to denote functions such that
$$|\pa_{t,x}^a r_m(t,x)|\le C_a \delta \<r\>^{-\rho-m-|a|}\ .$$
Then the equation $\Box_\gm u=F$ is equivalent to
\beeq\label{eq-LinWave}\Box u+ (g^{\al\be}-m^{\al\be})\pa_\al\pa_\be u=F-r_1 \pa u\equiv G\ .\eneq
With $f_j=r/(R+r)$ and $R\ge 1$,
through a direct calculation, it is easy to check from \eqref{520-4} that $Q_0$ is comparable to
\beeq\label{620-1}\frac{R}{(R+r)^{2}}\left[|\pa_r u|^2+
|\pt u|^2+\frac{R+r}{R} |\ang u|^2
+\left(1+(n-3)\frac{R+r}{R}\right)\frac{|u|^2}{r(R+r)}\right]\ .\eneq
By restricting the integral region to $\{r\le R\}$, we get
$$
\int_{S_T \cap \{r\le R\}} \frac{1}{R}\left(|\pa u|^2+\frac{|u|^2}{rR}\right) dx dt\les \int_{S_T} Q_0 dx dt\ .
$$
Together with Hardy's inequality and a cutoff argument\footnote{The author learned this simple but clever argument from Jason Metcalfe.} for the part $\{r<1\}$, i.e.,
$$\|u/r\|_{L^2_x(r< 1)}\le
\|\phi(r)u/r\|_{L^2_x}\les
\|\nabla(\phi(r)u)\|_{L^2_x}\les
\|u\|_{L^2_x(1<r<2)}+
\|\pa u\|_{L^2_x(r<2)},
$$ for some smooth cutoff function $\phi$,
 we see that
$$\|u\|_{LE}^2\les \sup_{R\ge 1}\int_{S_T} Q_0 dx dt\ .$$
Recalling \eqref{H1.1} and \eqref{eq-KSS3},  we could control the error terms $\int (Q-Q_0) dx dt$ by $\de\|u\|_{{LE}}^2$ (independent of $R\ge 1$), which yields
$$\|u\|_{LE}^2\les \sup_{R\ge 1}\int_{S_T} Q dx dt\ ,$$
for small enough $\de$.

Applying Lemma \ref{thm-KSS-0} and the standard energy inequality to the equation \eqref{eq-LinWave}, we could now obtain
\beeq\|u\|_{{LE}}\les \|\pa u(0)\|_{L^2_x}+
\|G\|_{L^1_t L^2_x+LE^*}
\ .\label{eq-612-3}
\eneq
Since $\rho>0$, we know that
$$\|r_1\pa u\|_{LE^*}\les \de \|u\|_{LE}\ .$$
This completes the proof of \eqref{eq-LE-g1-620}, in view of \eqref{eq-612-3} and smallness of $\de$.

\subsubsection{Nontrapping asymptotically Euclidean manifolds}

\begin{lem}\label{thm-LE-g2-620}
Let $n\ge 3$ and consider the linear problem $\Box_\gm u=F$ on the manifold $(\R^{1+n} , \gm)$ satisfying \eqref{H2} with $\rho=\min(\rho_1,\rho_2-1)>0$. Then for any positive constant $\mu>0$, we have the following local energy estimates,
\beeq\label{eq-LE-g2-620}
\|u\|_{LE}\les_\mu \|\pa u(0)\|_{L^2_x}+\|F\|_{L^1_t L^2_x +l^{1/2+\mu}_2 L^2_t L^2_x}\ .
\eneq
\end{lem}
\begin{prf}
It is essentially proved in Sogge and the author \cite{SoWa10}, where the authors used a cutoff argument to prove KSS estimates based on a weaker version of the local energy estimates of Bony and H\"afner \cite{BoHa}, which states that
for any $\mu>0$, \beeq\label{eq-LE-g2-620-2}
\|\pa u\|_{l^{-1/2-\mu}_2 L^2_t L^2_x}
+\|u\|_{l^{-3/2-\mu}_2 L^2_t L^2_x}
\les_\mu \|\pa u(0)\|_{L^2_x}+\|F\|_{L^1_t L^2_x +l^{1/2+\mu}_2 L^2_t L^2_x}\ .
\eneq
We need only to adapt that proof and notice that we could apply Lemma \ref{thm-LE-g1-620} for the part near spatial infinity. We leave the details to the reader. Here, it may be interesting to point out that, such an argument could yield \eqref{eq-LE-g2-620} for any $\mu>0$ from \eqref{eq-LE-g2-620-2} with certain $\mu_0>0$.
\end{prf}
\begin{rem}
As is clear from the proof, we do not need the radial assumption for $g_1$, in Lemmas \ref{thm-LE-g1-620}, \ref{thm-LE-g2-620} and \ref{thm-KSS-HY-inh}. The radial assumption will be required, however, for the higher order estimates in Lemma \ref{thm-LE-g3-high-620}, \ref{thm-KSS-HY-4}.
\end{rem}
\subsubsection{KSS estimates}
To prove almost global or long time existence, one may exploit the KSS estimates (see e.g. \cite{KSS02, MeSo06}). It is known that the KSS estimates are essentially contained in the local energy estimates (for $|x|\le 2+T$) and energy estimates (for $|x|\ge 2+T$), for which we record as the following lemma. One may consult  \cite{MeSo06} for the proof.
\begin{lem}
  \label{thm-KSS}
For any $\mu_1>0$, $\mu_2\in [0,1/2)$, there are positive constants $C_{\mu_1}$, $C_{\mu_2}$ and $C$, independent of $T>0$, such that
\beeq\label{eq-KSS3}
\|\pa u\|_{l^{-1/2-\mu_1}_2(L^2_T L^2_x)}+\| r^{-1}u\|_{l^{-1/2-\mu_1}_2(L^2_T L^2_x)}
\le C_{\mu_1} \|u\|_{LE ([0,T]\times\R^n)}\ ,
\eneq
\beeq\label{eq-KSS}
\|\pa u\|_{l^{-1/2}_2(L^2_T L^2_x)}+\|r^{-1} u\|_{l^{-1/2}_2(L^2_T L^2_x)}
\le C (\ln(2+T))^{1/2} \|u\|_{LE ([0,T]\times\R^n)}\ ,
\eneq
\beeq\label{eq-KSS2}
\|\pa u\|_{l^{-\mu_2}_2(L^2_T L^2_x)}+\| r^{-1} u\|_{l^{-\mu_2}_2(L^2_T L^2_x)}
\le C_{\mu_2} (2+T)^{1/2-\mu}\|u\|_{LE ([0,T]\times\R^n)}\ .
\eneq
\end{lem}

\subsubsection{higher order estimates}
\begin{lem}\label{thm-LE-g3-high-620}
Let $n\ge 3$ and consider the linear problem $\Box_\gm u=F$ on the manifold $(\R^{1+n} , \gm)$ satisfying \eqref{H1} or \eqref{H2}, with $\rho=\min(\rho_1,\rho_2-1)>0$. Then for small enough $\de$(in the case of \eqref{H1}), and for any positive constant $\mu>0$, we have the following higher order local energy estimates,
\beeq\label{eq-LE-g3-high-620}
\|u\|_{LE_k}\les \sum_{|a|\le k} \|\pa Y^a u(0)\|_{L^2_x}+\|Y^a F\|_{L^1_t L^2_x+l^{1/2+\mu}_2L^2_t L^2_x}\ .
\eneq
\end{lem}
\begin{prf}
The estimates for the asymptotically Euclidean manifolds are essentially contained in \cite{BoHa, SoWa10}, based on Lemma \ref{thm-LE-g2-620}. Here, we give only the proof for the asymptotically flat manifolds \eqref{H1}, in a manner similar to that of \cite{BoHa, SoWa10}, by induction. In this case, let us prove a stronger result
\beeq\label{eq-LE-g3-high-620-2}
\|u\|_{LE_k}\les \sum_{|a|\le k} \|\pa Y^a u(0)\|_{L^2_x}+\|Y^a F\|_{L^1_t L^2_x+LE^*}\ .
\eneq
At first,
Lemma \ref{thm-LE-g1-620} tells us that it is true for $k=0$. Assuming it is true for some $m\ge 0$, then
since $\Box_\gm u=F$ and \eqref{H1}, we have
\beeq\label{eq-relation1-2}\Box_\gm Y u=[\Box_\gm, Y]u+Y F=r_1 \pa^2 u+r_2 \pa u+Y F=r_1 \pa\nabla u+r_2 \pa u+Y F+r_1 F\ ,\eneq
where we have used the facts that $\pt^2 u=\Delta u+r_0 \pa \nabla u+F+r_0 F$ and the radial part of the metric has no effect in the commutator with rotational vector fields $\Omega$.
Applying \eqref{eq-LE-g3-high-620-2} to $Y u$ with $k=m$ and noticing \eqref{eq-relation1-2}, we see that
$\|Y u\|_{LE_m}$ can be controlled by
\begin{eqnarray*}
&\les&
 \sum_{|a|\le m} (\|\pa Y^a Y u(0,\cdot)\|_{L^2_x}+\|Y^a \Box_\gm Y u\|_{L^1_tL^2_x+LE^*})\\
&\les&
 \sum_{|a|\le m} (\|\pa Y^a Y u(0,\cdot)\|_{L^2_x}+\|Y^a ( r_1 \pa\nabla u+r_2 \pa u+Y F+r_1 F)\|_{L^1_tL^2_x+LE^*})\\
&\les&
  \sum_{|a|\le m+1} (\|\pa Y^a u(0,\cdot)\|_{L^2_x}+\|Y^a F\|_{L^1_tL^2_x+LE^*})+\sum_{|a|\le k-1}
 \|Y^a ( r_1 \pa\nabla u+r_2 \pa u)\|_{LE^*}\ .    \end{eqnarray*}
Recall that
$$Y^a (r_1 \pa\nabla u)=r_1\pa Y^a \nabla u+\sum_{|b|\le |a|}r_1\pa Y^b u\ ,$$
$$Y^a (r_2 \pa u)=r_2\pa Y^a  u+\sum_{|b|\le |a|}r_{2}\pa Y^b u\ ,$$
since $\rho>0$, it is easy to see that
$$\sum_{|a|\le m}
 \|Y^a ( r_1 \pa\nabla u+r_2 \pa u)\|_{LE^*}\les \sum_{|a|\le m+1}\|r_1\pa Y^a u\|_{LE^*}\les \de\sum_{|a|\le m+1}\|Y^a u\|_{LE}\ ,$$
 which completes the proof of \eqref{eq-LE-g3-high-620-2}, for small enough $\de$.
\end{prf}

\subsection{Local energy estimates, version 2}
In this subsection, for asymptotically flat manifolds, we will prove another version of the local energy estimates, which gives better control on the local part.
\begin{lem}[Local energy estimates]\label{thm-KSS-HY-inh}
Let $n\ge 3$ and consider the linear problem $\Box_\gm u=F$ on the manifold $(\R^{1+n} , \gm)$ satisfying \eqref{H1} with $\rho=\min(\rho_1,\rho_2-1)>0$. Then for any $\mu \in (0, 1/2)$, $\mu'-\mu\in (0, \rho/2]$,  and small enough $\de>0$,
\beeq\label{eq-KSS-HY-inh}
\|u\|_{\widetilde{LE}} \les \|\pa  u(0,\cdot)\|_{L^2_x}+\|\Box_\gm u\|_{L^1_t L^2_x + \widetilde{LE}^*}\ .
\eneq
\end{lem}

\begin{prf}
This result has essentially been proved in
Lemma 2.2 of \cite{HWY10}, using multipliers with $f_j=r/(2^j+r)$ ($j\ge 1$) and $f_0=(r/(1+r))^{2\mu}$ ($\mu\in (0,1/2)$).
Here, we use a single multiplier\footnote{The idea of constructing one multiplier is inspired by the discussion with Hans Lindblad when the author was working on \cite{LMSTW}.}  to give an alternative proof of the required local energy estimates.

The proof proceeds as that of Lemma \ref{thm-LE-g1-620} and we need only to point out the differences here.
For $\mu\in (0,1/2)$, $\mu'>\mu$, we set
\beeq\label{612-4}
f(r)=\left(\frac{r}{3+r}\right)^{2 \mu} \left(1-(3+r)^{-2 (\mu'-\mu)}\right)\ .
\eneq
Through a direct calculation, it is easy to check from \eqref{520-4} that $Q_0$ is comparable to
\beeq\label{620-2}r^{2\mu-1}\<r\>^{-2\mu'}\left[|\pa_r u|^2+
|\pt u|^2+\<r\>^{2(\mu'-\mu)} |\ang u|^2
+(1+(n-3)\<r\>^{2(\mu'-\mu)})r^{-2}u^2\right]\ ,\eneq
which controls the integrand in the local energy norm.
Applying Lemma \ref{thm-KSS-0} and the standard energy inequality to the equation \eqref{eq-LinWave}, we obtain
\beeq\|u\|_{\widetilde{LE}}\les \|\pa u(0)\|_{L^2_x}+
\|G\|_{L^1_t L^2_x+\widetilde{LE}^*}
\ .\label{612-3}
\eneq
Since $2\mu'\le \rho+2\mu$, we have
$$\|r_1\pa u\|_{\widetilde{LE}^*}\les \de \|u\|_{\widetilde{LE}}\ ,$$
and this completes the proof of \eqref{eq-KSS-HY-inh}, for small enough $\de$.
\end{prf}

To prove the radial Glassey conjecture, we need to obtain higher order local energy estimates involving the vector fields.
\begin{lem}\label{thm-KSS-HY-4}
Under the same assumption of Lemma \ref{thm-KSS-HY-inh}, we have the higher order local energy estimates for solutions to $\Box_\gm u=F$ with order $m\ge 0$:
\beeq\label{eq-KSS-HY-4'}
\| u\|_{\widetilde{LE}_m} \le C \sum_{|a|\le m} (\|\pa \nabla^a  u(0,\cdot)\|_{L^2_x}+\|\nabla^a F\|_{ L^1_tL^2_x+\widetilde{LE}^*})\ .
\eneq
In addition, we have
\beeq\label{eq-KSS-HY-310}
\sum_{|a|\le m}\| \pa^a u\|_{\widetilde{LE}} \le C \sum_{|a|\le m} (\|\pa \pa^a  u(0,\cdot)\|_{L^2_x}+\|\pa^a F\|_{ L^1_tL^2_x+\widetilde{LE}^*})\ .
\eneq
\end{lem}
The proof is basically the same as that of Lemma \ref{thm-LE-g3-high-620} and we omit the details here.

\subsubsection{Higher order local energy estimates}
Here, we record an alternative higher order local energy estimates, which shall be sufficient for us to adapt the current proof to obtain the same set of Theorems in \cite{HWY11}, for $n\ge 3$.
\begin{lem}\label{thm-KSS-HY-2}
Under the same assumption of Lemma \ref{thm-KSS-HY-inh}, we have the higher order local energy estimates for solutions to $\Box_\gm u=F$:
\beeq\label{eq-KSS-HY-2}
\|\nabla u\|_{\widetilde{LE}} \les\|\pa \nabla  u(0,\cdot)\|_{L^2_x}+\|\nabla F\|_{\widetilde{LE}^*}+\|r_1 F\|_{\widetilde{LE}^*}+\|r_2\pt u\|_{\widetilde{LE}^*}\ .
\eneq
In particular, letting $\|u\|_{\widehat{E}_2}=\|\nabla u\|_{E}+\de \|u\|_E$ and $\widehat{LE}_2$, $\widehat{LE}^*_2$ by similar fashion,
we have
\beeq\label{eq-KSS-HY-3}
\| u\|_{\widehat{LE}_2} \les\|\pa \nabla  u(0,\cdot)\|_{L^2_x}+\de \|\pa u(0,\cdot)\|_{L^2_x}+\|F\|_{\widehat{LE}_2^*}\ .
\eneq
\end{lem}
\begin{prf}
Since $\Box_\gm u=F$ and \eqref{H1}, we have
\beeq\label{eq-relation}\Box_\gm \nabla u=[\Box_\gm,\nabla ]u+\nabla F=r_1 \pa^2 u+r_2 \pa u+\nabla F=r_1 \pa\nabla u+r_2 \pa u+\nabla F+r_1 F\ ,\eneq
where we have used the fact that $\pt^2 u=\Delta u+r_0 \pa \nabla u+F+r_0 F$.
Notice that
$$\|r_1 \pa\nabla u+r_2 \nabla u\|_{\widetilde{LE}^*}\le C\de \|\nabla u\|_{\widetilde{LE}}\ ,$$
we see that, by applying Lemma \ref{thm-KSS-HY-inh},
\begin{eqnarray*}
     \|\nabla u\|_{\widetilde{LE}}
&\les& \|\pa \nabla  u(0,\cdot)\|_{L^2_x}+\|\Box_\gm \nabla u\|_{\widetilde{LE}^*}\\
          &\les &\|\pa \nabla  u(0,\cdot)\|_{L^2_x}+\de\|\nabla u\|_{\widetilde{LE}}+\|\nabla F\|_{\widetilde{LE}^*}+\|r_1 F\|_{\widetilde{LE}^*}+\|r_2\pa_t u\|_{\widetilde{LE}^*}\ ,
\end{eqnarray*} which gives us \eqref{eq-KSS-HY-2}.
For \eqref{eq-KSS-HY-3}, since $\|r_2 \pt u\|_{\widetilde{LE}^*}\le C\de \|u\|_{\widetilde{LE}}$, we need only to combine \eqref{eq-KSS-HY-inh} with \eqref{eq-KSS-HY-2}.
\end{prf}

\section{Glassey conjecture with dimension $3$}\label{sec-Gla-3}
In this section, we will prove Theorem \ref{thm-3}, mainly based on
Lemmas \ref{thm-Sobo} and \ref{thm-LE-g3-high-620}, for given $f$ and $g$ such that
\eqref{a-smallness} is satisfied.

As usual, we shall use iteration to give the proof. We set $u_0\equiv 0$ and recursively define $u_{k+1}$ be the solution to the linear equation
$$\Box_\gm u_{k+1}=F_p(u_k), u(0,x)=f(x), \pt u(0,x)=g(x).
$$

{\em Boundedness:}
By the smallness condition \eqref{a-smallness} on the data,  it follows from Lemma
\ref{thm-LE-g3-high-620} that there is a universal constant $C_1$ so that
\beeq\label{eq-620-15}\|u_1\|_{LE_2} \le C_1\ep, \
\|u_{k+1}\|_{LE_2} \le C_1\varepsilon + C_1 \sum_{|a|\le 2}\|Y^a F_p(u_k)\|_{L^1_tL^2_x}.\eneq

We shall argue inductively to prove that
\beeq\label{eq-unif2}\|u_{k+1}\|_{LE_2} \le 2C_1 \varepsilon.\eneq
By the above, it suffices to show
\beeq\label{eq-unif}
\sum_{|a|\le 2}\|Y^a F_p(u)\|_{L^1_tL^2_x}
\le \ep\ ,\eneq
for any $u$ with $\|u\|_{LE_2}\le 2 C_1 \ep\le 1$.

Notice that there exist smooth functions $b_i$, $1\le i\le 5$, such that
\begin{eqnarray*}
  Y^{\le 2} F_p(u) &= &b_1(u)|\pa u|^{p-1} Y^{\le 2} \pa u+b_2(u) |\pa u|^{p-2} (Y^{\le 1} \pa u)^2\\
&&+b_3(u) |\pa u|^{p-1} Y u Y^{\le 1} \pa u+b_4 (u) |\pa u|^p Y u Y u+b_5(u) |\pa u|^p Y^2 u\ .
  \end{eqnarray*}

By Lemma \ref{thm-Sobo},
\beeq\label{eq-bd} |\pa u|\les \frac{\|u\|_{E_{2}}}{\<r\>} ,\ |Y u|\les \|u\|_{E_{2}}\ .\eneq
Moreover, since $u(0)\in L^2_x$ and $\pt u\in L^\infty_t L^2_x$, we have $u(t)\in L^2_x$ for any $t$, which ensures that
\beeq\label{eq-bd-21}\|u(t)\|_{L^\infty_x}\le C\|\nabla u\|_{L^2_x}+C\|\nabla^2 u\|_{L^2} \le C \|u\|_{E_{2}}\ .\eneq Here, the constant $C$ is independent of the $L^2$ norm of $u(t)$.  By the boundedness of $u$, smoothness of $b_i$ and \eqref{eq-bd}, we see that
$$
|  Y^{\le 2} F_p(u) |\les |\pa u|^{p-1} (|Y^{\le 2} \pa u|+|Y^{\le 2} u|/\<r\>)+
|\pa u|^{p-2} |Y^{\le 1} \pa u|^2\ .$$

The first term can be dealt with as follows, by \eqref{eq-bd}, Lemma \ref{thm-Sobo} and Lemma \ref{thm-KSS} with $\mu_1=(p-2)/2$,
\begin{eqnarray*}
  &&\||\pa u|^{p-1} (|Y^{\le 2} \pa u|+|Y^{\le 2} u|/\<r\>)\|_{L^1_t L^2_x} \\
  &\les& \|\<r\>\pa u\|_{L_t^\infty L^\infty_x}^{p-2}\|\<r\>^{(3-p)/2}\pa u\|_{L^2_t L^\infty_x}\|  \<r\>^{-(p-1)/2}\left(|Y^{\le 2} \pa u|+\frac{|Y^{\le 2} u|}{\<r\>}\right)\|_{L^2_t L^2_x}\\
  &\les & \|u\|_{LE_{2}}^{p} \ .
\end{eqnarray*}  Similarly,
for the second term,
we get
\begin{eqnarray*}
\||\pa u|^{p-2} |Y^{\le 1} \pa u|^2\|_{L^1_t L^2_x}
  &\les&\|\<r\>\pa u\|_{L_t^\infty L^\infty_x}^{p-2} \|  \<r\>^{-(p-2)/2} Y^{\le 1} \pa u \|_{L^2_t L^4_x}^2\\
    &\les& \|u\|_{E_{2}}^{p-2}\|  \<r\>^{-(p-2)/2-1/2}Y^{\le 2} \pa u\|_{L^2_t L^2_{x}}^2\\
        &\les& \| u\|_{LE_{2}}^p\ .
\end{eqnarray*}
In conclusion, we see that there exists a constant $C_2$ such that
\beeq\label{0626-bd}\sum_{|a|\le 2}\|Y^a F_p(u)\|_{L^1_tL^2_x}\le C_2 \|u\|_{LE_2}^p\le C_2 (2C_1 \ep)^p\le \ep\eneq
for $\ep\le \ep_0$ with
$$C_2 (2C_1)^p \ep_0^p\le 1\ .$$
This finishes the proof of  \eqref{eq-unif} and so is the uniform boundedness \eqref{eq-unif2}.

{\em Convergence of the sequence $\{u_k\}$:}
Notice that
$$F_p(u)-F_p(v)=\sum_{\al} a_\al(u)(|\pa u|^p-|\pa v|^p)+|\pa v|^p (a_\al(u)-a_\al(v))\ .$$
Let $\ep_0>0$ be such that $2 C_4\ep_0\le 1$.
Since $a_\al$ are smooth functions, it is clear from \eqref{eq-bd-21}, \eqref{eq-unif2} that
$$|a_\al(u)-a_\al(v)|\les |u-v|,\ |a_\al(u)|\les 1\ , \forall u, v\in\{u_k\}.$$
Based on this observation, \eqref{eq-bd} and the fact $p>2$, together with Lemma \ref{thm-Sobo} and Lemma \ref{thm-KSS} with $\mu_1=(p-2)/2$, we see that
\begin{eqnarray*}
  &&\|F_p(u)-F_p(v)\|_{L^1_t L^2_x}\\
   &\les& \| (|\pa u|+|\pa v|)^{p-1}(|\pa (u- v)|+|u-v|/\<r\>)\|_{L^1_t L^2_x}\\
    &\les& (\|\<r\>\pa u\|_{L_t^\infty L^\infty_x}+\|\<r\>\pa v\|_{L_t^\infty L^\infty_x})^{p-2}\\
&&\times     (  \|\<r\>^{1/2-(p-2)/2}\pa u\|_{L^2_t L^\infty_x}+  \|\<r\>^{1/2-(p-2)/2}\pa v\|_{L^2_t L^\infty_x})\| u-v\|_{LE}\\
  &\les & (\|u\|_{LE_{2}}+\|v\|_{LE_2})^{p-1}\|u-v\|_{LE}\ ,
\end{eqnarray*}
that is, there exists a constant $C_3$ such that for any $u, v\in\{u_k\}$,
$$\|F_p(u)-F_p(v)\|_{L^1_tL^2_x}\le C_3 (\|u\|_{LE_2}
+\|v\|_{LE_2})^{p-1}\|u-v\|_{LE}
\le \|u-v\|_{LE}/(2C_1)\ ,$$
$$\|u_{k+1}-u_k\|_{LE}\le C_1\|F_p(u_k)-F_p(u_{k-1})\|_{L^1_tL^2_x}\le
\frac{1}2\|u_k-u_{k-1}\|_{LE}$$
provided that $\ep\le \ep_0$ with
$$2 C_1 C_3 (4C_1)^{p-1} \ep_0^{p-1}\le 1\ .$$

Together with the uniform boundedness \eqref{eq-unif2}, we find an unique global solution $u\in L^\infty_t H^3\cap Lip_t H^2$ with $\|u\|_{LE_2}\le 2 C_1 \ep$.

{\em Regularity of the solution:}
To complete the proof, it remains to prove the regularity of the solution $u\in C_t H^3\cap C^1_t H^2$. Since
$L^\infty_t H^3\cap Lip_t H^2\subset C_t H^2\cap C_t^1 H^1$, it suffices to prove $\pa u\in C_t H^2$.
For any $t_1>t_0\ge 0$, we set $w$ be the solution to the homogeneous wave equation
$$
\Box_\gm w=0, w(t_0)=u(t_0), \pt w(t_0)=\pt u(t_0)\ ,
$$
and so
$$
\Box_\gm (u-w)=F_p(u), (u-w)(t_0)=0, \pt (u-w)(t_0)=0\ ,
$$

 It is well-known that $w\in C_t H^3\cap C_t^1 H^2$.
Now,  by Lemma \ref{thm-LE-g3-high-620},
\begin{eqnarray*}
  \|u(t_1)-u(t_0)\|_{E_2}
&\le& \|u(t_1)-w(t_1)\|_{E_2}
+\|w(t_1)-u(t_0)\|_{E_2}\\
&\le&
\sum_{|a|\le 2}\|Y^a F_p(u)\|_{L^1_t([t_0,t_1]; L^2_x)}+o(t_1-t_0)\ .
\end{eqnarray*}
Since $Y^a F_p(u)\in L^1_t  L^2_x$ (see \eqref{0626-bd}), it is clear that $\|Y^a F_p(u)\|_{L^1_t([t_0,t_1]; L^2_x)}=o(t_1-t_0)$ and so is $\|u(t_1)-u(t_0)\|_{E_2}$, which completes the proof of the regularity.

\section{High dimensional wave equations}\label{sec-Gla-4}
In this section, for high dimensional wave equations, we will prove Theorem \ref{thm-1}, based on
Lemmas \ref{thm-Sobo} and \ref{thm-LE-g3-high-620}, for given $f$ and $g$ such that
\eqref{a-smallness2} is satisfied.

We will proceed as in Section \ref{sec-Gla-3} to give the proof and we give only the proof of the uniform boundedness of the iteration series $u_k$.

Based on \eqref{eq-620-15},
we shall argue inductively to prove that
\beeq\label{eq-unif2-2}\|u_{k+1}\|_{LE_m} \le 2C_1 \varepsilon ,\eneq
for which it suffices to show
\beeq\label{eq-unif-2}
\sum_{|a|\le m}\|Y^a F_p(u)\|_{L^1_T L^2_x}
\le \ep\ ,\eneq
for any $u$ with $\|u\|_{LE_m}\le 2 C_1 \ep\le 1$.

We are assuming $a_\al (u)=a_\al$,  and so
$|Y^{\le m} F_p(u)|$ is controlled by terms of the following form
\beeq\label{eq-typical-1} |\pa u|^{p-j}\Pi_{i=1}^j |Y^{a_i} \pa u| \ , 1\le a_i\le a_{i+1},\ \sum_{i=1}^j a_i\le m\ .\eneq
Here, we have abused the notation $a_i$ with the order $|a_i|$.
We will consider two cases separately: 1) $j\le 1$, 2) $j\ge 2$.

Case 1), $j\le 1$. In this case,
by  Lemma \ref{thm-Sobo}, we see that
\begin{eqnarray*}
\||\pa u|^{p-1}Y^{\le m} \pa u\|_{L^1_t L^2_x}&\le&
\|\<r\>^{(n-1)/2}\pa u\|_{L^\infty_t L^\infty_x}^{p-2}
\|\<r\>^{(n-1)/2-(n-1)(p-1)/4}\pa u\|_{L^2_t L^\infty_x}\\
&&\times\|\<r\>^{-(n-1)(p-1)/4} Y^{\le m} \pa u\|_{L^2_t L^2_x}
\\
&\les&
\|Y^{\le m}\pa u\|_{L^\infty_t L^2_x}^{p-2}
\|\<r\>^{-(n-1)(p-1)/4}Y^{\le m}\pa u\|_{L^2_t L^2_x}^2\\
\end{eqnarray*}

Case 2) $j\ge 2$. Since $a_{1}\ge 1$ and $\sum a_i\le m$, we get $a_i\in [1,[n/2]]$ and
$$c_i=m-a_i\in [1,[n/2]]\ .$$
Fix $\ep_i\in (0,c_i)$, such that \beeq\label{eq-typical-2}\sum \ep_i=-\frac n2(j-1)+\sum c_i \ .\eneq
Notice that it is possible since
$$
-\frac n2(j-1)+\sum c_i=
-\frac n2(j-1)+mj-\sum a_i\ge (j-1)(m-\frac n2)>0\ .$$
We can define $p_i\in (2,\infty)$ through
$$n(\frac 12-\frac 1{p_i})=c_i-\ep_i\in (0,n/2)\ .$$
Here, we notice that \eqref{eq-typical-2} ensures that $\sum 1/p_i=1/2$.

Now, we  set
$d_i=(n-1)(p-j)/(2j)+(n-1)(1/2-1/p_i)$ for $i< j$,
and $d_j=-\sum_{i=1}^{j-1} d_i$. Then
$$d_j-(n-1)(p-j)/(2j)-(n-1)(1/2-1/p_j)=
-(n-1)(p-1)/2\ .
$$
Applying Lemma \ref{thm-Sobo} and H\"older's inequality, we could control $\||\pa u|^{p-j} \Pi_{i=1}^j Y^{a_i} \pa u\|_{L^1_t L^2_x}$ as follows
\begin{eqnarray*}
&\le&
\|\<r\>^{(n-1)/2}\pa u\|_{L^\infty_t L^\infty_x}^{p-j}
\Pi_{i=1}^{j-2}\|\<r\>^{d_i-(n-1)(p-j)/(2j)}Y^{a_i}\pa u\|_{L^\infty_t L^{p_i}_x}\\
&&\times
\|\<r\>^{d_{j-1}-(n-1)(p-j)/(2j)-(n-1)(p-1)/4}Y^{a_{j-1}}\pa u\|_{L^2_t L^{p_{j-1}}_x}
\\&&\times
\|\<r\>^{d_j-(n-1)(p-j)/(2j)+(n-1)(p-1)/4} Y^{a_j} \pa u\|_{L^2_t L^{p_j}_x}
\\
&\les&
\|Y^{\le m}\pa u\|_{L^\infty_t L^2_x}^{p-j}
\Pi_{i=1}^{j-2}\|\<r\>^{(n-1)(1/2-1/p_i)}Y^{a_i}\pa u\|_{L^\infty_t L^{p_i}_x}\\
&&\times
\|\<r\>^{(n-1)(1/2-1/p_{j-1})-(n-1)(p-1)/4}Y^{a_{j-1}}\pa u\|_{L^2_t L^{p_{j-1}}_x}
\\&&\times
\|\<r\>^{(n-1)(1/2-1/p_j)-(n-1)(p-1)/4} Y^{a_j} \pa u\|_{L^2_t L^{p_j}_x}
\\
&\les&
\|Y^{\le m}\pa u\|_{L^\infty_t L^2_x}^{p-j}
\Pi_{i=1}^{j-2}\|Y^{\le m}\pa u\|_{L^\infty_t L^{2}_x}\\
&&\times
\|\<r\>^{-(n-1)(p-1)/4} Y^{\le m} \pa u\|_{L^2_t L^{2}_x}^2
\\
&\les&
\|Y^{\le m}\pa u\|_{L^\infty_t L^2_x}^{p-2}
\|\<r\>^{-(n-1)(p-1)/4}Y^{\le m}\pa u\|_{L^2_t L^2_x}^2
\end{eqnarray*}

In summary, we have proved that
$$\|Y^{\le m} F_p(u)\|_{L^1_t L^2_x}
\les\|Y^{\le m}\pa u\|_{L^\infty_t L^2_x}^{p-2}
\|\<r\>^{-(n-1)(p-1)/4}Y^{\le m}\pa u\|_{L^2_t L^2_x}^2
.$$
When $p>1+2/(n-1)$, by
Lemma \ref{thm-KSS} with $\mu_1=(n-1)(p-1)/4-1/2>0$, it is controlled by $\|u\|_{LE_m}^p$ and so is \eqref{eq-unif-2},  with small enough $\ep$.
For the remaining situation, that is $n=3$ and $p=2$, we have $(n-1)(p-2)/4=1/2$, and we get from \eqref{eq-KSS} in Lemma \ref{thm-KSS}  that
$$\|Y^{\le m} F_p(u)\|_{L^1_T L^2_x}
\les
\ln(2+T) \|u\|_{LE_m}^2\ ,$$
which yields \eqref{eq-unif-2}, for $T\le exp(c/\ep)$ with small enough $\ep$ and $c$.

\section{Radial Glassey conjecture}
In this section, we give the proof for Theorem \ref{thm-sub}, by using
Lemma \ref{thm-trace} and Lemma \ref{thm-KSS-HY-4}.

We set $u_0\equiv 0$ and recursively define $u_{k+1}$
to be the solution to the linear equation
\beeq\label{3iterate}
  \Box_\gm u_{k+1} = c_1|\pt u_k|^p+c_2|\nabla u_k|^p\equiv G_p(u_{k}),    u(0,x)= f,
  \pt u(0,x) = g.
\eneq

As in Section \ref{sec-Gla-4}, we give only the proof of the uniform boundedness of the iteration series $u_k$ with respect to $\widetilde{LE}_1$, which could be reduced to the proof of
\beeq\label{eq-radG-1}\|u\|_{\widetilde{LE}_1}\le 2C_1\ep \Rightarrow \|G_p(u)\|_{\widetilde{LE}_1^*}\le \ep
\eneq
by Lemma \ref{thm-KSS-HY-4} with $m=1$.

Recall from \eqref{eq-trace-rad} in Lemma \ref{thm-trace} with $s=1$,  we have for radial functions $u$ (see e.g.  \cite{HWY11} Lemma 2.4),
\beeq\label{eq-trace-rad2}
\|r^{(n-2)/2}\<r\>^{1/2}  \pa u\|_{L^\infty_x}
\les \|r^{(n-2)/2}\<r\>^{1/2}  \pa u\|_{L^\infty_r L^2_\omega}
\les  \|\nabla^{\le 1}\pa u\|_{ L^2}\ .
\eneq
So
\begin{eqnarray*}
\|G_p(u)\|_{\widetilde{LE}_1^*}&\les&
\| r^{1/2-\mu}\<r\>^{\mu'} \nabla^{\le 1} G_p(u)\|_{L^2 L^2}
\\
&\les&\| r^{1/2-\mu}\<r\>^{\mu'} |\pa u|^{p-1} \nabla^{\le 1} \pa u\|_{L^2 L^2}\\
&\les&
\|r^{(n-2)/2}\<r\>^{1/2} \pa u\|_{L^\infty_{t,x}}^{p-1}
\| r^{1/2-\mu-(n-2)(p-1)/2}\<r\>^{\mu'-(p-1)/2}  \nabla^{\le 1} \pa u\|_{L^2 L^2}\\
&\les&
\| u\|_{LE_1}^p
\end{eqnarray*}
provided that
$1/2-\mu-(n-2)(p-1)/2\ge -1/2+\mu$,
$1/2-\mu-(n-2)(p-1)/2+\mu'-(p-1)/2\le -1/2+\mu-\mu'$,
which is true if we set
$$\mu=\frac{1}2- \frac{n-2}4 (p-1)\ ,\ \mu'=\frac{p-1}4\ .$$
Then by choosing $\ep>0$ small enough, we get \eqref{eq-radG-1}, and so is Theorem \ref{thm-sub}.

\section{Global existence for radial data}
In this section, we give the proof for Theorem \ref{thm-add}, filling some gap between Theorem \ref{thm-1} and \ref{thm-sub}, by using
Lemma \ref{thm-trace}, Lemma \ref{thm-KSS-HY-4}, as well as Lemma \ref{thm-Sobo}, \ref{thm-wHLS}.

\subsection{Dimension $n=5$}\label{sec-7.1}
At first, let us give the proof for $n=5$ and $p\in [2,3)$. We give only the proof of the uniform boundedness here, which, by \eqref{eq-KSS-HY-310}, amounts to the proof of
\beeq\label{eq-radGh-310} \|\nabla^{\le 2} \pa u(0)\|_{L^2_x}\le \ep
\Rightarrow
\|\pa^{\le 2} \pa u(0)\|_{L^2_x}\les \ep\ ,
\eneq
\beeq\label{eq-radGh-1}\|\pa^{\le 2} u\|_{\widetilde{LE}}\les \ep \Rightarrow \|\pa^{\le 2} G_p(u)\|_{\widetilde{LE}^*+L^1_t L^2_x}\le \ep\ ,
\eneq
for radial functions and small enough $\ep$.

For \eqref{eq-radGh-310}, we observe from the equation that 
$$|\pa^{\le 2} \pa u(0)|\les
|\nabla^{\le 2} \pa u(0)|+
|\pa u(0)|^{p-1}|\nabla^{\le 1} \pa u(0)|+
|\pa u(0)|^{2p-1}\ .
$$
Since $p\in [2,3)$, by Sobolev embeddings $H^1\subset L^2\cap L^{10/3}$ and $H^2\subset L^2\cap L^{10}$, we see that
$$
\|\pa^{\le 2} \pa u(0)\|_{L^2}\les
\|\nabla^{\le 2} \pa u(0)\|_{L^2}+
\|\nabla^{\le 2} \pa u(0)\|_{L^2}^p+
\|\nabla^{\le 2} \pa u(0)\|_{L^2}^{2p-1}\ ,
$$ which completes the proof of \eqref{eq-radGh-310}.

Recall from Lemma \ref{thm-Sobo} with $k=1$, that 
$$\|\<r\>^{(n-1)/n}u\|_{L^{2n/(n-2)}_x}\les \|Y^{\le 1} u\|_{L^2_x},
$$
from which, together with \eqref{eq-trace-rad} with $s=1$, and H\"older's inequality, yields
\beeq\label{eq-radGh-2}\|r^{1/4}\<r\>^{3/4}u\|_{L^{4}_r L^3_\omega}\les \|Y^{\le 1} u\|_{L^2_x}, n=5.
\eneq
Here, we introduced the mixed-norm $L^\rho_r L^s_\omega$ with respect to the polar coordinates $x=r\omega$, $\omega\in \Sp^{n-1}$
$$
\|f\|_{L^\rho_r L^s_\omega}=\left(\int\|f(r \cdot)\|_{L^s_\omega}^\rho r^{n-1}dr\right)^{1/\rho}
,
$$
with obvious modification for $\rho=\infty$.
A cutoff argument  gives us that for any $m_1, m_2\in \R$, we have
\beeq\label{eq-radGh-2'}\|r^{1/4+m_1}\<r\>^{3/4+m_2}u\|_{L^{4}_r L^3_\omega}\les \|r^{m_1}\<r\>^{m_2}Y^{\le 1} u\|_{L^2}
+\|r^{m_1-1}\<r\>^{m_2}u\|_{L^2}
, n=5.
\eneq
Actually, introducing a radial cutoff function $\phi\in C_0^\infty(\R^n)$ with $\phi=1$ for $|x|\le 2$ and $\phi=0$ for $|x|\ge 3$, and applying \eqref{eq-radGh-2}, we get
\begin{eqnarray*}
\|r^{1/4+m_1}\<r\>^{3/4+m_2} \phi u
\|_{L^4}
&\les&
\|r^{1/4}\<r\>^{3/4}(\phi r^{m_1}u)\|_{L^4}
\\&\les&
\sum_{|a|\le 1} \| Y^a(r^{m_1}\phi u)\|_{L^2( r\le 3)}
\\&\les&
\|r^{m_1-1}   u\|_{L^2(r\le 3)}+ \| r^{m_1} Y u\|_{L^2( r\le 3)}
\ ,
\end{eqnarray*}
and for $j\ge 1$,
\begin{eqnarray*}
\|r^{1/4+m_1}\<r\>^{3/4+m_2} u
\|_{L^4(j+1\le r\le j+2)}
&\les&
j^{m_1+m_2}\|
r^{1/4}\<r\>^{3/4}u\|_{L^4(j+1\le r\le j+2)}
\\&\les&
j^{m_1+m_2}\sum_{|a|\le 1} \| Y^a u\|_{L^2(j\le r\le j+3)}
\ .
\end{eqnarray*}
The above two inequalities together imply \eqref{eq-radGh-2'} if we  $l^4$-sum over $j\ge 0$.

Also, by Lemma \ref{thm-trace} with $s=2$,  we have for radial functions $u$ (see e.g.  \cite{HWY11} Lemma 2.4),
\beeq\label{eq-radGh-3}
\|r^{1/2}\<r\>^{3/2} \pa u\|_{L^\infty_x}\les
\|r^{1/2}\<r\>^{3/2} \pa u\|_{L^\infty_r L^2_\omega}
\les \|\nabla^{\le 2} \pa u\|_{ L^2}, n=5\ .
\eneq

There are two kinds of typical terms in $\pa^{\le 2} G_p(u)$, i.e.,
$|\pa u|^{p-1} \pa^{\le 2}\pa u$ and $|\pa u|^{p-2} |\pa^2 u|^2$.

The first term can be dealt with as follows, by \eqref{eq-radGh-3},
\begin{eqnarray*}
  &&\||\pa u|^{p-1}\pa^{\le 2} \pa u\|_{\widetilde{LE}^*} \\
  &\les&
  \|r^{1/2}\<r\>^{3/2}\pa u\|_{L^\infty_t L^\infty_x}^{p-1}
  \|r^{1/2-\mu}\<r\>^{\mu'}(r^{1/2}\<r\>^{3/2})^{-(p-1)}\pa^{\le 2} \pa u\|_{L_t^2 L^2_x}\\
  &\les&  \|\nabla^{\le 2}\pa u\|_{L^\infty_x L^2_x}^{p-1}  \|r^{-1/2+\mu}\<r\>^{-\mu'}\pa^{\le 2} \pa u\|_{L_t^2 L^2_x}\\
&\les & \|\pa^{\le 2} u\|_{\widetilde{LE}}^{p} \ ,
\end{eqnarray*}
if \beeq\label{eq-radGh-4}
\mu= (3-p)/4, \mu'-\mu= p-3/2\ .
\eneq
Similarly,
for the second term, by \eqref{eq-radGh-3}, we get
\begin{eqnarray*}
&&\||\pa u|^{p-2} |\pa^2 u|^2\|_{L^1_t L^2_x} \\
  &\les&\|r^{1/2}\<r\>^{3/2}\pa u\|_{L_t^\infty L^\infty_x}^{p-2} \|
  r^{-(p-2)/4}\<r\>^{-3(p-2)/4}\pa^2 u \|_{L^2_t L^4_x}^2\\
    &\les& \|\nabla^{\le 2}\pa u\|_{L^\infty_t L^2_x}^{p-2} \|
  r^{-(p-2)/4}\<r\>^{-3(p-2)/4}\pa^2 u \|_{L^2_t L^4_x}^2\ .
\end{eqnarray*}
Here, for the last term, we distinguish two cases: $\pa^2 u=\pa \pt u$ and $\pa^2 u=\nabla^2 u$. For the first case, since $u$ is radial, $\pt u$ is also radial and $\nabla \pt u=\frac{x}r\pa_r \pt u $. By \eqref{eq-radGh-4} and \eqref{eq-radGh-2'}, we have
\begin{eqnarray*}
&&\|
  r^{-(p-2)/4}\<r\>^{-3(p-2)/4}\pa \pt u \|_{L^2_t L^4_x}\\
&\le &\|
  r^{-(p-2)/4}\<r\>^{-3(p-2)/4}\pa_{t,r} \pt u \|_{L^2_t L^4_r L^3_\omega}\\
&\le&
\|  r^{-(p-1)/4}\<r\>^{-3(p-1)/4}\nabla^{\le 1} \pa_{t,r}\pt u\|_{L^2_t L^2_{x}}+
    \|  r^{-(p+3)/4}\<r\>^{-3(p-1)/4}\pa_{t,r}\pt u\|_{L^2_t L^2_{x}}\\
&\le&
\|  r^{-(p-1)/4}\<r\>^{-3(p-1)/4}\pa^{\le 2}\pt u\|_{L^2_t L^2_{x}}+
    \|  r^{-(p+3)/4}\<r\>^{-3(p-1)/4}\pa\pt u\|_{L^2_t L^2_{x}}\\
    &\les&
    \|\pa^{\le 2} u\|_{\widetilde{LE}}\ .
\end{eqnarray*}

For the second case, it is a little more involved, due to the fact that $\nabla u$ is nonradial.
Anyway, in this case, we observe that $|\nabla^2 u|\les |\pa^2_r u|+|\pa_r u|/r$. As a consequence, 
\begin{eqnarray*}
&&\|
  r^{-(p-2)/4}\<r\>^{-3(p-2)/4}\nabla \nabla u \|_{L^4_x}\\
&\les &
\|  r^{-(p-2)/4}\<r\>^{-3(p-2)/4}\pa_r^2  u \|_{L^4_r L^3_\omega}
+\|  r^{-(p-2)/4-1}\<r\>^{-3(p-2)/4}\pa_r  u \|_{L^4_r L^3_\omega}
\\
&\les&
\|  r^{-(p-1)/4}\<r\>^{-3(p-1)/4}\nabla^{\le 1} \pa_r^2 u\|_{ L^2_{x}}+
    \|  r^{-(p-1)/4-1}\<r\>^{-3(p-1)/4}\pa_r^2 u\|_{ L^2_{x}}\\
&&+
\|  r^{-(p-1)/4-1}\<r\>^{-3(p-1)/4} \nabla^{\le 1}\pa_r  u \|_{L^2_x}
+
\|  r^{-(p-1)/4-2}\<r\>^{-3(p-1)/4}\pa_r  u \|_{L^2_x}
\\
&\les&
\|  r^{-(p-1)/4}\<r\>^{-3(p-1)/4}\nabla^{\le 3} u\|_{ L^2_{x}}+
    \|  r^{-(p-1)/4-1}\<r\>^{-3(p-1)/4}\nabla^2 u\|_{ L^2_{x}}\\
&&+
\|  r^{-(p-1)/4-2}\<r\>^{-3(p-1)/4}\nabla  u \|_{L^2_x}
\\
&\les&
\|  r^{-(p-1)/4}\<r\>^{-3(p-1)/4}\nabla^{\le 3} u\|_{ L^2_{x}}+
    \|  r^{-(p-1)/4-1}\<r\>^{-3(p-1)/4}\nabla^2 u\|_{ L^2_{x}}\\
&&+
\|  r^{-(p-1)/4}\Delta(\<r\>^{-3(p-1)/4}\nabla  u )\|_{L^2_x}
\\&\les&
\|  r^{-(p-1)/4}\<r\>^{-3(p-1)/4}\nabla^{\le 3} u\|_{ L^2_{x}}+
    \|  r^{-(p-1)/4-1}\<r\>^{-3(p-1)/4}\nabla^2 u\|_{ L^2_{x}}\ ,
\end{eqnarray*}
where, in the second to the last inequality, we have used Lemma \ref{thm-wHLS} with $k=2$.
\footnote{The argument used here is inspired by an unpublished work of the author, in collaboration with Kunio Hidano and Kazuyoshi Yokoyama.
}

In conclusion, we have proved $$\||\pa u|^{p-1}\pa^{\le 2} \pa u\|_{\widetilde{LE}^*} +\||\pa u|^{p-2} |\pa \pa u|^2\|_{L^1_t L^2_x} \les 
\|\pa^{\le 2} u\|_{\widetilde{LE}}^p\ ,$$ which completes the proof of \eqref{eq-radGh-1}.

\subsection{Even Spatial Dimension $n\ge 4$}
In this subsection, we turn to the proof for the even spatial dimension $n\ge 4$. As before, we give only the proof of
\beeq\label{eq-radGh4-310} \|\nabla^{\le n/2} \pa u(0)\|_{L^2_x}\le \ep, \|\pa u(0)\|_{L^\infty_x}\les 1
\Rightarrow
\|\pa^{\le n/2} \pa u(0)\|_{L^2_x}\les \ep\ ,
\eneq
\beeq\label{eq-radGh4-1}\|\pa^{\le n/2}u\|_{\widetilde{LE}}\les\ep \Rightarrow \|\pa^{\le {n/2}} G_p(u)\|
_{\widetilde{LE}^*}\le \ep\ , 
\eneq
for radial functions and small enough $\ep$.

For \eqref{eq-radGh4-310}, by equation and a tedious calculation, we can find that
$$|\pa^a \pa u(0)|\les |\pa u(0)|^{|a|(p-1)+1}+
\sum_{|a_i|\ge 1, \sum |a_i|\le |a|} |\pa u(0)|^{(|a|-\sum |a_i|)(p-1)-m+1}\Pi_{i=1}^m |\nabla^{a_i} \pa u(0)|\ .$$
The $L^2$ norm of the first term on the right is trivial. For the remaining terms, since $\pa u(0)$ is bounded, we need only to control the $L^2$ norm of  $\Pi_{i=1}^m |\nabla^{a_i} \pa u(0)|$. For this purpose, we introduce $q_i\in [2,\infty)$, satisfying 
$$\sum \frac{1}{q_i}=\frac 12, \frac{n}{q_i}\ge |a_i|,$$
which is possible since $\sum |a_i|\le |a|\le n/2$. 
Then by Sobolev embeddings $H^{n/2-|a_i|}\subset L^{q_i}$ and H\"{o}lder's inequality, we get
\begin{eqnarray*}
\|\Pi_{i=1}^m |\nabla^{a_i} \pa u(0)|\|_{L^2}& \les & \Pi_{i=1}^m \|\nabla^{a_i} \pa u(0)\|_{L^{q_i}}\\
& \les & \Pi_{i=1}^m \|\nabla^{a_i} \pa u(0)\|_{H^{n/2-|a_i|}}\\
 & \les & \|\pa^{\le n/2} \pa u(0)\|_{L^2_x}^m \ ,
 \end{eqnarray*}
 which yields the desired result \eqref{eq-radGh4-310}.

By Lemma \ref{thm-trace}, 
we have for radial functions,
\beeq\label{eq-radGh4-3}
\|r^{\epsilon}\<r\>^{(n-1)/2-\epsilon} \pa u\|_{L^\infty_x}\les
\|r^{\epsilon}\<r\>^{(n-1)/2-\epsilon} \pa u\|_{L^\infty_r L^2_\omega}\les \|\pa u\|_{ H^{n/2}}, \forall \epsilon>0\ .
\eneq

Since $G_p(u)=c_1 |\pt u|^p+c_2|\nabla u|^p$, $|\pa^{\le {n/2}} G_p(u)|$ is controlled by terms of the following form
\beeq\label{eq-typical4-1} |\pa u|^{p-j}\Pi_{i=1}^j |\pa^{a_i} \pa u| \ , 1\le a_i\le a_{i+1},\ \sum_{i=1}^j a_i\le {n/2}\ .\eneq
We will consider two cases separately: 1) $j\le 1$, 2) $j\ge 2$.

Case 1), $j\le 1$. In this case,
with $\mu=\frac{1}4$, $\mu'=\frac{(n-1)(p-1)-1}4$, we set
$\vep=1/(2(p-1))$. By
\eqref{eq-radGh4-3},
, we have
\begin{eqnarray*}
&&\||\pa u|^{p-1}\pa^{\le {n/2}} \pa u\|_{\widetilde{LE}^*}\\
&\le&
\|r^\vep\<r\>^{(n-1)/2-\vep}\pa u\|_{L^\infty_t L^\infty_x}^{p-1}
\|r^{1/2-\mu-\vep (p-1)}\<r\>^{\mu'-[(n-1)/2-\vep](p-1)}
\pa^{\le {n/2}} \pa u\|_{L^2_t L^2_x}
\\
&\les&
\|\pa^{\le n/2}\pa u\|_{L^\infty_t L^2_x}^{p-1}
\|\pa^{\le n/2}u\|_{\widetilde{LE}}\ .
\end{eqnarray*}

Case 2) $j\ge 2$. Since $a_{1}\ge 1$ and $\sum a_i\le {n/2}$, we get $a_i\in [1,n/2-1]$ and
$$c_i=n/2-a_i\in [1,n/2-1]\in (0,n/2)\ .$$

Fix $\ep_i\in [0,c_i)$, such that \beeq\label{eq-typical4-2}\sum \ep_i=-\frac n2(j-1)+\sum c_i \ .\eneq
Notice that it is possible since
$$
-\frac n2(j-1)+\sum c_i=
-\frac n2(j-1)+j \frac n 2-\sum a_i\ge 0\ .$$
We can define $p_i\in (2,\infty)$ through
$$n\left(\frac 12-\frac 1{p_i}\right)=c_i-\ep_i\in \left(0,\frac n2\right)\ .$$
Here, we notice that \eqref{eq-typical4-2} ensures that $\sum 1/p_i=1/2$.
Setting $\vep =1/(2(p-j))$,
we have $$\frac{n-1-2\vep}2(p-j)+\sum_{i=1}^{j}(n-1)(\frac 12-\frac{1}{p_i})=2\mu'\ .$$
Then, applying Lemma \ref{thm-Sobo}, \eqref{eq-radGh4-3} and H\"older's inequality, we get
\begin{eqnarray*}
&&\||\pt u|^{p-j} \Pi_{i=1}^j \pa^{a_i} \pa u\|_{\widetilde{LE}^*}\\
&\le&
\|r^\vep\<r\>^{(n-1)/2-\vep}\pa u\|_{L^\infty_t L^\infty_x}^{p-j}
\Pi_{i=1}^{j-1}\|
\<r\>^{(n-1)(1/2-1/p_i)}\pa^{a_i}\pa u\|_{L^\infty_t L^{p_i}_x}
\\&&\times
\|r^{1/2-\mu-\vep (p-j)}
\<r\>^{\mu'-(n-1-2\vep)(p-j)/2
-\sum_{i=1}^{j-1}
(n-1)(1/2-1/p_i)
} \pa^{a_j} \pa u\|_{L^2_t L^{p_j}_x}
\\
&\les&
\|\nabla^{\le n/2}\pa u\|_{L^\infty_t L^2_x}^{p-j}
\|r^{-1/2+\mu}
\<r\>^{-\mu'+(n-1)(1/2-1/p_j)}
 \pa^{ a_j} \pa u\|_{L^2_t L^{p_j}_x}\\
 &&\times
\Pi_{i=1}^{j-1}\|Y^{\le c_i}\pa^{ {a_i}}\pa u\|_{L^\infty_t L^2_x}
\\
&\les&
\|\pa^{\le n/2}\pa u\|_{L^\infty_t L^2_x}^{p-1}
\|r^{-1/2+\mu}
\<r\>^{-\mu'+(n-1)(1/2-1/p_j)}
 \pa^{ a_j} \pa u\|_{L^2_t L^{p_j}_x}\ ,
\end{eqnarray*}
where in the last inequality, we have used the fact that $|Y^a \pa u|\les |\nabla^{\le a} \pa u|$ for radial functions. It remains to prove
\beeq\label{claim}\|r^{-1/2+\mu}
\<r\>^{-\mu'+(n-1)(1/2-1/p_j)}
 \pa^{ a_j} \pa u\|_{L^2_t L^{p_j}_x}\les \|\pa^{\le n/2}u\|_{\widetilde{LE}}\ .\eneq
In fact, by Lemma \ref{thm-Sobo} and \ref{thm-wHLS}, we know that
\begin{eqnarray*}
&&\|r^{-1/2+\mu}
\<r\>^{-\mu'+(n-1)(1/2-1/p_j)}
 \pa^{ a_j} \pa u\|_{L^{p_j}_x}\\
 & \les &
\|Y^{\le c_j}(r^{-1/2+\mu}
\<r\>^{-\mu'}
 \pa^{ a_j} \pa u)\|_{L^{2}_x} 
  \\
 & \les & \|r^{-1/2+\mu}
\<r\>^{-\mu'}
 \pa^{\le n/2} \pa u\|_{L^{2}_x} 
+\sum_{k=1}^{c_j}\sum_{|b|\le c_j-k} \|r^{-1/2+\mu-k}
\<r\>^{-\mu'}
 \pa^{\le a_j+b} \pa u\|_{L^{2}_x} \\
 & \les & \|r^{-1/2+\mu}
\<r\>^{-\mu'}
 \pa^{\le n/2} \pa u\|_{L^{2}_x} 
 +\|r^{-3/2+\mu}
\<r\>^{-\mu'}
 \pa^{\le n/2-1} \pa u\|_{L^{2}_x} 
 \\
 &&
+\sum_{3\le 2l+1\le c_j}\sum_{|b|\le c_j-2l} \|r^{-1/2+\mu-2l-1}
\<r\>^{-\mu'}
 \pa^{\le a_j+b} \pa u\|_{L^{2}_x} \\
&&+\sum_{2\le 2l\le c_j}\sum_{|b|\le c_j-2l} \|r^{-1/2+\mu-2l}
\<r\>^{-\mu'}
 \pa^{\le a_j+b} \pa u\|_{L^{2}_x} \\
  & \les & \|r^{-1/2+\mu}
\<r\>^{-\mu'}
 \pa^{\le n/2} \pa u\|_{L^{2}_x} 
 +\|r^{-3/2+\mu}
\<r\>^{-\mu'}
 \pa^{\le n/2}  u\|_{L^{2}_x} 
 \\
 &&
+\sum_{3\le 2l+1\le c_j}\sum_{|b|\le c_j-2l} \|r^{-1/2+\mu-1}
\Delta^l (\<r\>^{-\mu'}
 \pa^{\le a_j+b} \pa u)\|_{L^{2}_x} \\
&&+\sum_{2\le 2l\le c_j}\sum_{|b|\le c_j-2l} \|r^{-1/2+\mu}
\Delta^l(\<r\>^{-\mu'}
 \pa^{\le a_j+b} \pa u)\|_{L^{2}_x} \\
   & \les & \|r^{-1/2+\mu}
\<r\>^{-\mu'}
 \pa^{\le n/2} \pa u\|_{L^{2}_x} 
 +\|r^{-3/2+\mu}
\<r\>^{-\mu'}
 \pa^{\le n/2}  u\|_{L^{2}_x} \ .
 \end{eqnarray*}

In summary, we have proved that
$$\|\pa^{\le n/2} G_p(u)\|_{\widetilde{LE}^*}
\les
\|\pa^{\le n/2}\pa u\|_{L^\infty_t L^2_x}^{p-1}
\|\pa^{\le n/2}u\|_{\widetilde{LE}}
$$
which is sufficient to prove
\eqref{eq-radGh4-1},
  with small enough $\ep$.

\medskip 
{\bf Acknowledgements.}
The author would like to thank the anonymous referee for valuable comments, which helped improve this manuscript.


\end{document}